\documentclass[10pt]{amsart}
\usepackage{amsmath,amstext,amscd, amssymb,txfonts, epsfig, psfrag}

\theoremstyle{plain}
 \newtheorem{MainThm}{Theorem}
 \newtheorem{MainCor}[MainThm]{Corollary}

\newtheorem{thm}{Theorem}[section]
\newtheorem{lemma}[thm]{Lemma}

\newtheorem{prop}[thm]{Proposition}

\theoremstyle{definition}
\newtheorem{defn}[thm]{Definition}

\newtheorem{Remark}[thm]{Remark}

\newtheorem{Open questions}[thm]{Open questions}
\newtheorem{Open question}[thm]{Open question}
\newtheorem{Open problems}[thm]{Open problems}
\newtheorem{Open problem}[thm]{Open problem}

\def\cal{\mathcal}

\def\Bbb{\mathbb}
\def\bar{\overline}

\def\Z{\Bbb{Z}}
\def\N{\Bbb{N}}

\def\Reals{\Bbb{R}}

\def\ni{\noindent}

\def\Area{\hbox{\rm Area}}
\def\Cay{\hbox{\it Cay}}

\def\IDiam{\hbox{\rm Diam}}

\def\IDiam{\hbox{\rm IDiam}}

\def\FL{\hbox{\rm FL}}
\def\FFL{\hbox{\rm FFL}}
\def\FFFL{\hbox{\rm FFFL}}
\def\F+L{\hbox{$\textup{F}\!_+\textup{L}$}}

\def\ssm{\smallsetminus}

\def\ms{\medskip}

\def\onto{{\kern3pt\to\kern-8pt\to\kern3pt}}

\def\<{\langle}
\def\>{\rangle}
\def\|{{\ |\ }}

 \def\AA{\cal A}

 \def\RR{\cal R}

 \def\PP{\cal P}
 
\def\QQ{\cal Q}

\newcommand{\set}[1]{\left\{#1\right\}}
\newcommand{\restricted}[1]{\left|_{#1} \right.}
\newcommand{\abs}[1]{\left|#1\right|}
\renewcommand{\ni}{\noindent}
\renewcommand{\ss}{\smallskip}
\renewcommand{\ms}{\medskip}
\newcommand{\bs}{\bigskip}
\def\*{^{\star}}

\begin{document}

\markboth{M.R.Bridson and T.R.Riley}{Free and fragmenting filling length}

\title{Free and Fragmenting Filling Length}

\author{M.\ R.\ Bridson
 and T.\ R.\ Riley}

\date{December 8, 2005; revised May 16, 2006.  To appear in the Journal of Algebra}

\begin{abstract}
\ni The filling length of an edge-circuit $\eta$ in the Cayley
2-complex of a finitely presented group is the least integer $L$ such that there is a combinatorial null-homotopy of $\eta$ down
to a basepoint through loops of length at most $L$. We introduce similar notions in which the null-homotopy is not required to fix a basepoint, and in which the contracting loop is allowed to bifurcate.  We exhibit groups in which the resulting filling invariants exhibit dramatically different behaviour to the standard notion of filling length.  We also define the corresponding filling invariants for Riemannian manifolds and translate our results to this setting.   \ss
\\
\footnotesize{\ni \textbf{2000 Mathematics Subject
Classification:  20F65, 20F10, 53C23}  \\ \ni \emph{Key words and phrases:} filling
length, null-homotopy, van~Kampen diagram}
\end{abstract}

\maketitle

\section{Introduction} \label{intro}

\ni Consider a  vertical cylinder $C \subseteq \mathbb{R}^3$  of height $h$ 
whose base has  diameter $d \ll h$. Let $S$ be the surface formed by the curved portion of $C$ and the disc capping off its top.  Topologically, $S$ is a closed 2-disc.  The loop $\partial S$ can be homotoped in $S$ to a constant
loop through loops of length at most $\pi d$ by lifting it up the
cylinder and then contracting it across the top of $C$.
However, if we insist on keeping a basepoint on $\partial S$ fixed in the course of the null-homotopy then we will encounter far longer loops, some of length at least $2h$.

In this article we will bring to light similar contrasts between basepoint-fixed and basepoint-free null-homotopies for loops in the Cayley 2-complex $Cay^2(\PP)$ of a finite presentation $\PP$ of a group $\Gamma$. Words $w$ that represent $1$ in
$\Gamma$ (\emph{null-homotopic} words) correspond to edge-circuits $\eta_w$ in $\Cay^2(\PP)$.
The filling length $\FL(w)$ of $w$ was defined by Gromov \cite{Gromov} and in a combinatorial context is the minimal length $L$ such that there is a basepoint-preserving
 \emph{combinatorial null-homotopy}  of $\eta_w$ through loops of length at most $L$.  (A closely related notion called $\textup{LNCH}$  was considered by Gersten in \cite{Gersten5}.)  
We define $\FFL(w)$, the \emph{free filling length} of $w$, likewise but without holding a
basepoint fixed, and $\FFFL(w)$, the \emph{fragmenting free filling length} of $w$, by also allowing the contracting loops to bifurcate.  Detailed definitions are in Section~\ref{filling length}.  

We construct a finite presentation in which $\FL(w)$ and $\FFL(w)$ differ dramatically for an  infinite sequence of null-homotopic words of increasing length.  [Our conventions are $[a,b]= a^{-1}b^{-1}ab$, $a^b = b^{-1}ab$ and $a^{-b}= b^{-1}a^{-1}b$. For $f,g:\N \to \N$ or $f,g:\Reals \to \Reals$, we write $f \preceq g$ when there exists $C>0$ such that for all $l$ we have $f(l) \leq Cg(Cl+C)+Cl+C$, which gives an equivalence relation expressing qualitative agreement of the growths of $f$ and $g$: write $f \simeq g$ if and only if $f \preceq g$ and $g \preceq f$.  We compare functions $\N \to \N$ with functions  $\Reals \to \Reals$ by extending the former to $\Reals$ by considering them to be constant on the intervals $[n,n+1)$ for all $n \in \N$.] 

\begin{MainThm} \label{Thm 1}
Let $\Gamma$ be the group given by the aspherical presentation
$$\QQ \ := \ \langle \ a,b, r, s, t \ \mid \ a^ba^{-2}, [t,a], [r, at], [r,s], [s,t]  \ \rangle.$$
For $n \in \N$ define $w_n :=  [s,a^{-b^n}  r a^{b^n}].$ Then
$w_n$ is null-homotopic in $\QQ$, has length $\ell(w_n) = 8n+8$, and $\FFL(w_n) \simeq n$, but $\FL(w_n) \simeq  2^n$. 
\end{MainThm}

We assume the reader is familiar with \emph{van~Kampen diagrams} (\cite{Bridson6} is a recent survey); they can be thought of as \emph{combinatorial homotopy discs} for loops in the Cayley 2-complex of a presentation. We show that there are van~Kampen diagrams  $\hat{\Delta}_n$  for $w_n$ that, owing to geometry like that in our cylinder example, have $\FFL(\hat{\Delta}_n) \preceq n$ and $\FL(\hat{\Delta}_n) \simeq 2^n$.   Indeed, we will show that every van~Kampen diagram $\Delta$ for $w_n$ has intrinsic  diameter $ \succeq 2^n$; that is, there is a pair of vertices $x,y$ such that every edge-path in $\Delta^{(1)}$ connecting $x$ to $y$ has length  $ \succeq 2^n$.  It will then follow that $\FL(w_n) \succeq 2^n$.

The area $\Area(w)$ (resp.\ intrinsic diameter $\IDiam(w)$ \,) of $w$ is the 
least integer $K$ such that there is a van~Kampen diagram for $w$ with $K$ 2-cells (resp.\ with intrinsic diameter $K$).   For $\textup{M} =  \Area, \IDiam, \FL, \FFL$ or $\FFFL$ one defines \emph{filling functions} $\textup{M} : \N \to \N$ for finitely presented groups:  $$\textup{M}(n) \ := \ \max\set{ \ \textup{M}(w) \mid w \textup{ null-homotopic and } \ell(w) \leq n \  }.$$  (The argument of $\textup{M}$ --- diagram, word or
integer --- determines its meaning;  the potential for ambiguity is tolerated as it spares us from a terminology  over-load.)  In the case $\textup{M} = \Area$, the function $M(n)$ is the \emph{Dehn function}.

In spite of Theorem~\ref{Thm 1}, the filling functions $\FL$ and $\FFL$ for $\QQ$ are $\simeq$-equivalent: we will see that any van~Kampen diagram for $w'_n:=  [s, \  a^{-b^n} r a^{ b^n} r a^{b^n} r^{-1} a^{-b^n}]$ has two \emph{peaks} and the savings that can be made by escaping the basepoint are no longer significant.  On the other hand, if we allow our loops to bifurcate then they can pass over peaks independently, and so $\FFFL$ exhibits markedly different behaviour. 

\begin{MainThm} \label{Thm 2}
The filling functions $\IDiam, \FL, \FFL, \FFFL: \N \to \N$ for
 $\QQ$ satisfy 
\begin{eqnarray*}
\IDiam(n) \ \simeq \ \FL(n) \ \simeq \ \FFL(n) \ \simeq \ \Area(n) & \simeq & 2^n \\
\FFFL(n) & \simeq & n.
\end{eqnarray*}
\end{MainThm}

The relations $\IDiam(n) \simeq \Area(n)$ and $\FL(n) \simeq \Area(n)$ are noteworthy in the non-hyperbolic setting.  In contrast, for groups
 where $\Area(n) \simeq n^{\alpha}$ for some $\alpha \geq 2$ one
 knows that $\IDiam(n) \preceq n^{\alpha-1}$ -- see \cite{GR1}.  

\ms

Theorems~\ref{Thm 1} and \ref{Thm 2} are proved in Section~\ref{main proof}, modulo a number of auxiliary propositions postponed to Section~\ref{aux}.  The remaining sections are dedicated to establishing the credentials of $\FL$, $\FFL$ and $\FFFL$ for inclusion in the pantheon of filling invariants: we relate them to other filling functions. we interpret them in terms of algorithmic complexity, and we relate them to functions concerning null-homotopies of loops in metric spaces (particularly, Riemannian manifolds).

\ms

 $\FL$ can be thought of as a space-complexity measure in that $\FL(w)$ is the minimal $L$ such that $w$ can be converted to the empty word through a sequence of words of length at most $L$, each obtained from the previous by free reduction, free expansion, or applying a relator  (see \cite{GR1}).  We will show in Section~\ref{filling length} that  $\FFL$ and $\FFFL$ can be interpreted in a similar way: by allowing  an
additional operation of conjugation we get $\FFL(w)$, and further including the move that replaces a word $w=uv$ by a pair of words $u,v$ we get $\FFFL(w)$.  This point of view is useful for calculations and allows us to prove that, as for $\FL$, given a finite presentation, $\Area(n)$ is at most an exponential of $\FFFL$.  This and other relationships are spelt out in the following theorem, which shows, in particular, that $\QQ$ provides an example of $\Area(n)$ outgrowing $\FFFL(n)$ as extremely as is possible.  

\begin{MainThm} \label{exp bound}
Let $\PP$ be a finite presentation.   There is a constant $C$, depending only on $\PP$, such that for all $n \in \N$ the Dehn, filling length, free filling length, and fragmenting free filling length functions $\Area, \FL, \FFL, \FFFL: \N \to
\N$ of $\PP$ satisfy 
\begin{eqnarray*}
\Area(n) & \leq & C^{\mbox{$\FFFL(n)$}} \\ 
\FFFL(n) & \leq & \FFL(n) \ \leq \ \FL(n) \ \leq \ C\, \Area(n) + n.
\end{eqnarray*}
\end{MainThm}

In Section~\ref{qi} we prove that, up to $\simeq$-equivalence,
  $\FL(n)$,  $\FFL(n)$ and  $\FFFL(n)$ are all quasi-isometry invariants
  among finitely presented
  groups and so, in particular, do not depend on a choice of finite presentation of a group -- cf.\ \cite[Theorem~8.1]{GR4}.

\begin{MainThm} \label{qi thm}
If $\PP$ and $\PP'$ are finite presentations of quasi-isometric groups then 
 $$\FL_{\PP} \simeq \FL_{\PP'},  \ \  \ \ \FFL_{\PP} \simeq \FFL_{\PP'},  \ \textit{and } \ \   \FFFL_{\PP} \simeq \FFFL_{\PP'}.$$  
\end{MainThm}

In  Section~\ref{rc} we define analogous filling invariants $\FL_X(n)$,  $\FFL_X(n)$ and  $\FFFL_X(n)$ concerning null-homotopies of rectifiable loops in arbitrary metric spaces $X$, and we prove:

\begin{MainThm} \label{translate}
Suppose a group $\Gamma$ with finite presentation $\PP = \langle \AA \mid \RR \rangle$ acts properly and cocompactly by
isometries on a simply connected geodesic metric space $X$ for which there exist $\mu, L > 0$ such that every loop of length less than $\mu$ admits a based null-homotopy of filling length less than $L$. 
Then $\FL_{\PP} \simeq \FL_X$, $\FFL_{\PP} \simeq \FFL_X$ and $\FFFL_{\PP} \simeq \FFFL_X$.
\end{MainThm}

As the universal cover of any closed connected Riemiannian manifold satisfies these conditions and (as is well known) every finitely presentable group is the fundamental group of such a manifold, we can use Theorem~\ref{translate} and its analogues for $\Area$ and $\IDiam$ (proved in \cite{Bridson6, BR1}) to obtain from Theorems~\ref{Thm 1} and \ref{Thm 2}: 

\begin{MainCor} \label{cor}
There exists a closed connected Riemannian manifold $M$ such that 
\begin{eqnarray*}
\IDiam_{\widetilde{M}}(n) \ \simeq \ \FL_{\widetilde{M}}(n) \ \simeq \ \FFL_{\widetilde{M}}(n) \ \simeq \ \Area_{\widetilde{M}}(n) & \simeq & 2^n \\
\FFFL_{\widetilde{M}}(n) & \simeq & n.
\end{eqnarray*}
Moreover, there is an infinite sequence of loops $c_n$ in $\widetilde{M}$ such that $\ell(c_n) \to \infty$,  $\FFL_{\widetilde{M}}(c_n) \simeq n$, and $\FL_{\widetilde{M}}(c_n) \simeq 2^n$.
\end{MainCor}

\ni [Strictly speaking, the final part is not a direct corollary of the prior results, but it follows from the methods of  Section~\ref{rc}: the $w_n$ of Theorem~\ref{Thm 1} can be used to construct \emph{word-like} loops $c_n$ with $\ell(c_n) \simeq \ell(w_n)$; van~Kampen diagrams witnessing to the upper bounds $\FFL(w_n) \preceq n$ and $\FL(w_n) \preceq 2^n$ can be translated into fillings of $c_n$ showing  $\FFL_{\widetilde{M}}(c_n) \preceq n$ and $\FL_{\widetilde{M}}(c_n) \preceq 2^n$; and every null-homotopy of $c_n$ in $\widetilde{M}$ has filling length  $\succeq 2^n$ because the filling-disc corresponding to a null-homotopy facilitates a van~Kampen diagram for $w_n$ with $\simeq$ similar filling length.]

\ms

Natural questions that remain open include the following.

\begin{Open problem}
Does there exist a finite presentation for which $\FL(n) \nsimeq \FFL(n)$? 
\end{Open problem}

\begin{Open problem}
Does there exist a finite presentation for which $\IDiam(n) \nsimeq \FL(n)$? 
\end{Open problem}

In reference to the first of these problems we note that for a finite presentation $\langle \AA \mid \RR \rangle$ of a group $\Gamma$, the presentation $\langle \AA \, \cup \, \set{t} \mid \RR \rangle$ for $\Gamma \ast \Z$ satisfies $\FL(n) \simeq \FFL(n)$ for the following reason.  If $w_n$ is a word over $\PP$ with $\FL(w_n) = \FL(n)$ then $\FFL([w_n,t]) \geq 2+ \ell(w_n) + \FL(w_n)$ since any diagram for $[w_n,t]$ is a diagram for $w_n$ joined to a diagram for ${w_n}^{-1}$ by a $t$-edge, and the most efficient way to shell such a diagram (from the point-of-view of $\FFL$) is to shell the $w_n$ diagram, then collapse the $t$-edge, and then shell the ${w_n}^{-1}$ diagram.  

The second problem was posed by Gromov in \cite[\S5C]{Gromov}, where it is also attributed to Casson. A negative answer would imply that the \emph{double exponential upper bound} \cite{Cohen, Gersten2} on the Dehn (i.e.\ area) function in terms of $\IDiam(n)$ could be improved to a single exponential  
using the single exponential upper bound \cite{GR1, Gromov} for $\Area(n)$ in terms of $\FL(n)$.

\bs
\ni \emph{Acknowledgement.}  The second author is grateful for support from NSF grant  DMS--0404767 and for the hospitality of the Centre de Recerca Matem\`atica in Barcelona during the writing of this article.

\section{Three notions of filling length} \label{filling length}

\ni Let $\Delta$ be a \emph{diagram}; that is, a finite, planar,
contractible, combinatorial 2-complex; i.e.\ a van~Kampen diagram bereft of all group theoretic decorations. Before discussing filling length, 
we recall a combinatorial notion of null-homotopy 
  called a \emph{shelling}  from \cite{GR1}.

\begin{figure}[ht]
\psfrag{1}{1-cell}%
\psfrag{2}{2-cell}%
\psfrag{f}{fragmentation}%
\psfrag{c}{collapse}%
\psfrag{e}{expansion}%
\centerline{\epsfig{file=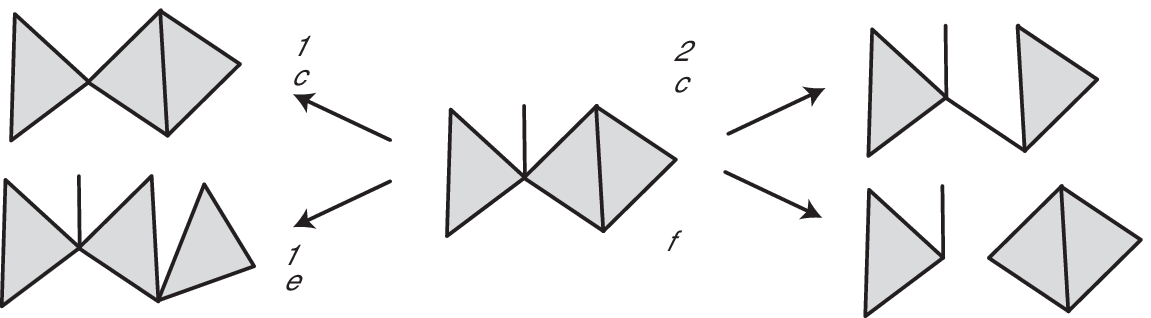}} \caption{Shelling
moves.} \label{shelling move fig}
\end{figure}

\begin{defn}\label{shelling defn}
A \emph{shelling} $\mathcal{S}= (\Delta_i)$ of $\Delta$ is a
sequence of diagrams
$$\Delta \ = \ \Delta_0, \, \Delta_1, \, \dots, \, \Delta_m,$$ in which each $\Delta_{i+1}$ is obtained from $\Delta_i$
by one of the \emph{shelling moves} defined below and depicted in
Figure~\ref{shelling move fig}.
\begin{itemize}
\setlength{\itemsep}{0pt} \setlength{\parsep}{0pt}

\item \emph{1-cell collapse.} Remove a pair $(e^1,e^0)$ where
$e^1$ is a 1-cell with $e^0 \in \partial e^1$ and $e^1$ is
attached to the rest of $\Delta_i$ only by one of its end vertices
$\neq e^0$.  (We call such an $e^1$ a \emph{spike}.)

\item \emph{1-cell expansion.} Cut along some 1-cell $e^1$ in $\Delta_i$ that has a vertex $e^0$ in $\partial
\Delta_i$, in such a way that $e^0$ and $e^1$ are doubled.

\item \emph{2-cell collapse.} Remove a pair $(e^2,e^1)$ where
$e^2$ is a 2-cell which has some edge $e^1 \in (\partial e^2 \cap
\partial \Delta_i)$.  The effect on the boundary circuit is to replace
$e^1$ with $\partial e^2 \ssm e^1$.
\end{itemize}

We say that the shelling $\mathcal{S}$ is \emph{full} when $\Delta_m$ is a single vertex.
A \emph{full shelling to a base vertex $\star=\Delta_m$} on $\partial \Delta$ is a full shelling in which $\star$ is preserved
throughout the sequence $(\Delta_i)$. In particular, in
every \emph{1-cell collapse} $e^0 \neq \star$, and in every
\emph{1-cell expansion} on $\Delta_i$ where $e_0=\star$ a choice is made as to
which of the two copies of $e_0$ is to be $\star$ in $\Delta_{i+1}$.

We define a \emph{full fragmenting shelling} $\mathcal{S}$ of
$\Delta$ by adapting the definition of a full shelling by allowing each $\Delta_i$ to be a disjoint union of finitely many diagrams, 
insisting that $\Delta_m$ be a set of vertices, and allowing one
extra type of move:

\begin{itemize}
\item \emph{Fragmentation.} $\Delta_{i+1}$ is the disjoint union
of $\Delta'_i$ and $\Delta''_i$, where $\Delta_i = \Delta'_i \cup \Delta''_i$ and $\Delta'_i \cap
\Delta''_i$ is a single vertex.
\end{itemize}
\end{defn}

For a shelling $\mathcal{S}$, define $\ell(\mathcal{S}):= \max_i
\ell(\partial \Delta_i)$, where $\ell(\partial \Delta_i)$ denotes the sum of the lengths
of the boundary circuits of the components of $\Delta_i$.  Then the \emph{filling
length} $\FL(\Delta, \star)$, the
\emph{free filling length} $\FFL(\Delta)$, and the
\emph{fragmented free filling length} $\FFFL(\Delta)$ are the
minima of $\ell(\mathcal{S})$ as $\mathcal{S}$ ranges over all full
shellings to $\star$, full shellings, and all full fragmenting free shellings
of $\Delta$, respectively.  This notation emphasises the fact that $\FL(\Delta,\star)$ is defined with respect to a base vertex $\star \in \partial \Delta$ but $\FFL(\Delta)$ and $\FFFL(\Delta)$ are not.

\ms

Let $\PP= \langle \AA \mid \RR \rangle$  be a finite presentation
of a group $\Gamma$.  Define the filling length $\FL(w)$ of a null-homotopic word $w$ in
a presentation $\PP= \langle \AA \mid \RR \rangle$ by $$\FL(w):=
\min \set{ \FL(\Delta,\star) \mid \Delta \textup{ a van~Kampen diagram
for } w},$$ and $\FFL(w)$ and $\FFFL(w)$ likewise.

\bs
A sequence $\mathcal{S}= w_0,
\ldots, w_m$ of null-homotopic words is a \emph{null-sequence} if each
$w_{i+1}$ is obtained from $w_i$ by one of the following moves:
\begin{itemize}
\setlength{\itemsep}{0pt} \setlength{\parsep}{0pt}

\item \emph{Free reduction.} Remove a subword $aa^{-1}$ or
$a^{-1}a$ from $w_i$, where $a \in \AA$.

\item \emph{Free expansion.} This is the inverse of a free reduction.

\item \emph{Application of a relator.} Replace a subword $u$ of
$w_i$ by a word $v$ such that a cyclic conjugate of $uv^{-1}$ is
in $\RR^{\pm 1}$.
\end{itemize}

\ni We define two more moves:

\begin{itemize}
\setlength{\itemsep}{0pt} \setlength{\parsep}{0pt}
\item \emph{Cyclic conjugation.} Replace $w_i$ by a cyclic
permutation.
\item \emph{Fragmentation.}  Replace a word $w=uv$ by a pair of
words $u$, $v$. (In effect, insert a letter into $w$
that represents a blank space.)
\end{itemize}
To employ the fragmentation move we must generalise our definition of a
\emph{null-sequence} so that each $w_i$ is a finite sequence of
words, and when we perform any of the operations listed above we execute it on one of the words in $w_i$.  

The proof of the following reassuring lemma is straightforward.  The ``only if'' part is well known (indeed, freely reducing, freely expanding, and applying relators suffices).  The ``if'' part can be proved by an easy induction on the number of fragmentation moves used. 

\begin{lemma}
A word $w$ over $\PP$ represents the identity if and only if it can be reduced to a sequence of empty words by free reductions, free expansions, applying relators, cyclically
conjugating, and fragmenting.
\end{lemma}

For a null-sequence $\mathcal{S}$, define $\ell(\mathcal{S}):= \max_i \ell(w_i)$ where, if fragmentation moves are employed, $\ell(w_i)$ is the sum of the lengths of words in the sequence $w_i$.  
Proposition~1 in \cite{GR1} says that for all null-homotopic words $w$, we have
$\FL(w) = \min_{\mathcal{S}} \ell(\mathcal{S})$, quantifying over null-sequences $\mathcal{S}$ for $w$ that employ  \emph{free reductions}, \emph{free expansions} and \emph{applications of relators}.  We add the following.

\begin{prop} \label{null seq prop}
Quantifying over all null-sequences $\mathcal{S}$ for a null-homotopic word $w$, where \emph{free reduction}, \emph{free expansion}, \emph{applications of relators}, and \emph{cyclic conjugation} are allowed, we have  
$\FFL(w)= \min_{\mathcal{S}} \ell(S).$
If, additionally, we allow \emph{fragmentations} we get $\FFFL(w)= \min_{\mathcal{S}} \ell(S).$
\end{prop}

\begin{proof}  
The proof in \cite{GR1} that $\min_{\mathcal{S}} \ell(\mathcal{S}) \leq \FL(w)$ is straightforward --- the 
words around the boundary of the van~Kampen diagrams in the course of a full shelling form a null-sequence.
Each word in the sequence of boundary words in the course of a \emph{free} shelling of a van~Kampen diagram is only defined up to cyclic conjugation, as a base vertex is not kept fixed during the shelling.  
It is then easy to see that $\min_{\mathcal{S}} \ell(S) \leq \FFL(w)$, quantifying over all all null-sequences $\mathcal{S}$ for $w$ that use \emph{free reduction}, \emph{free expansion}, \emph{applications of relators}, and \emph{cyclic conjugation}.  Introducing \emph{fragmentation} moves into the shelling, produces \emph{fragmentations} in the corresponding null-sequence, and we see  $\min_{\mathcal{S}} \ell(S) \leq \FFFL(w)$.     

The reverse bounds require more care.  Given a null-sequence $\mathcal{S}$ for $w$ involving \emph{applications of relations} and \emph{free expansions} and \emph{reductions}, we seek to construct a van~Kampen diagram with a shelling during which the lengths of the boundary circuit remain at most $\ell(\mathcal{S})$.  This can be done by starting with an edge-loop labelled by $w$, and filling it in by attaching a 2-cell on every \emph{application of a relator}, by attaching a 1-cell on every \emph{free expansion}, by folding together two adjacent 1-cells on every \emph{free reduction}.  However it is possible that the resulting complex will not be planar: 2-spheres or other cycles may be pinched off (for example, when an inverse pair $aa^{-1}$ is inserted and then removed).  Removing these cycles gives a van~Kampen diagram $\Delta$ with $\FL(\Delta) \leq \ell(w)$.  This is explained carefully in \cite{GR1}.  If $\mathcal{S}$ also uses \emph{cyclic conjugation} moves then no extra complications are added to the construction of $\Delta$.  If $\mathcal{S}$ also uses \emph{fragmentations} then the corresponding move in the course of the construction of $\Delta$ is to identify two vertices so that an inner boundary circuit is changed from a topological circle to a figure-eight.    \end{proof}

\ms
\begin{Remark} (\emph{Filling length and space complexity}.)
Envisaging a null-sequence to be the course of a calculation
on a Turing tape, we see that $\FL(w)$ is the non-deterministic
space complexity of the following approach to solving the word
problem for $\PP$: write $w$ on the tape and then exhaustively
apply relators and perform free reductions and free expansions.  A
sequence of moves that converts $w$ into the empty word amounts to a
proof that $w$ represents 1 in $\PP$ and $\FL(w)$ is the minimal
upper bound on the number of places on the tape that have to be
used in the calculation (see \cite{GHR} for more details). If we
allow cyclic conjugation then the non-deterministic space
complexity is $\FFL(w)$. If we also include fragmentation then the non-deterministic space complexity is $\FFFL(w)$ plus the maximum number of blank spaces separating the words; that is, between  $\FFFL(w)$ and $2\ \FFFL(w)$. 
\end{Remark}

\bs

The following inequalities are easy  
\begin{eqnarray}
\FFFL(w) \  \leq \  \FFL(w) &  \leq &  \FL(w) \ \leq K\, \Area(w) + \ell(w), \label{ineq1} \\
\IDiam(w) & \leq & \FL(w),
\end{eqnarray} where $K$ is a constant depending only on $\PP$.   
(See \cite{GR1} or \cite{Gromov}.)
Similar inequalities, for example $\IDiam(n) \leq \FL(n)$ for all $n$, relate the corresponding filling functions (defined in Section~\ref{intro}).

It is well known \cite{GR1, Gromov} that there is a
constant $C$, that depends only on $\PP$, such that the Dehn
function $\Area: \N \to \N$ of $\PP$ satisfies $$\Area(n) \leq
C^{\mbox{$\FL(n)$}}$$ for all $n$.  This is essentially the Time-Space
bound of algorithmic complexity: the number of different words of length
$\FL(n)$ on an alphabet of size $C$ is $C^{\mbox{$\FL(n)$}}$, and as $w$ admits a  null-sequence that does not include repeated words, $C^{\mbox{$\FL(n)$}}$
is an upper bound on the number of times relations are applied
in the null-sequence and hence on $\Area(w)$.  The
same proof applies to $\FFFL$ and so with (\ref{ineq1}) gives Theorem~\ref{exp bound}.

\section{Proofs of Theorems~\ref{Thm 1}~and~\ref{Thm 2}} \label{main proof}

\begin{proof}[Proof of Theorem~\ref{Thm 1}.]  
We exploit a technique due to the first author in \cite{Bridson} to show that the 
intrinsic diameter $\IDiam(w_n)$ of $w_n=[s, a^{- b^n} r a^{ b^n}]$ is at least $2^n$. 

\begin{figure}[ht]
\psfrag{D}{$\Delta$}
\psfrag{S}{$\Sigma$}
\psfrag{u0}{$u_0$}
\psfrag{mu}{$\mu$}
\psfrag{u}{$u$}
\psfrag{v}{$v$}
\psfrag{a}{$a$}
\psfrag{b^n}{$b^n$}
\psfrag{a2n}{$a^{2^n}$}
\psfrag{t2n}{$t^{2^n}$}
\psfrag{at}{$(at)^{2^n}$}
\psfrag{T}{$s$}
\psfrag{t}{$t$}
\psfrag{tau}{$r$}
\psfrag{a}{$a$}
\psfrag{b}{$b$}
\psfrag{T}{$s$}
\psfrag{t}{$t$}
\psfrag{tau}{$r$}
\psfrag{s}{$\star$}
\centerline{\epsfig{file=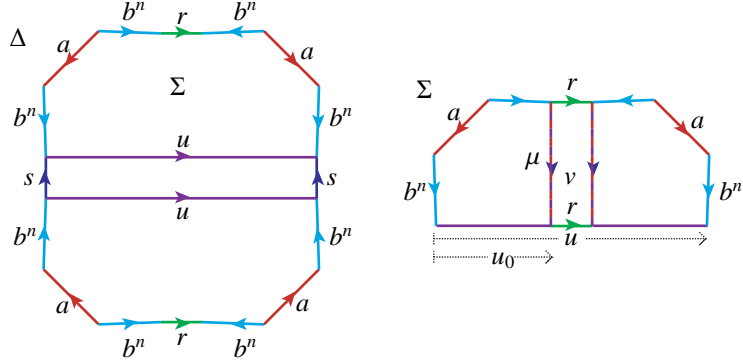}} \caption{A van~Kampen
diagram $\Delta$ for $w_n$ with a subdiagram $\Sigma$.} \label{Delta}
\end{figure}

Suppose $\pi: \Delta \to
\Cay^2(\QQ)$ is a van~Kampen diagram for $w_n$. Figure~\ref{Delta} is a schematic depiction of $\Delta$ and Figure~\ref{Delta_n} shows an explicit example when $n=3$.  We
seek an edge-path $\mu$ in $\Delta$, labelled by a word
in which the exponent sum of the letters $t$ is $2^n$. A $s$-corridor through $\Delta$ connects the two
letters $s$ in $w_n$.  Along each side of this corridor we read a word $u$ in $\set{r^{\pm 1}, t^{\pm 1}}^{\star}$.  Let $\Sigma$ be a subdiagram of $\Delta$ with boundary made up of one side of the $s$-corridor and a portion of $\partial \Delta$ labelled $a^{-b^n}r a^{ b^n}$. An $r$-corridor in $\Sigma$ joins the $r$ in $ a^{-   b^n} r  a^{b^n}$ to some edge-labelled
$r$ in $u$.  Let $u_0$ be the prefix of $u$ such that the
letter immediately following $u_0$ is this $r$, and then
let $\mu$ be the edge-path along the side of the $r$-corridor
running from the vertex at the end of $a^{-  b^n}$ to the
vertex at the end of $u_0$.  Let $v \in  \set{(at)^{\pm 1}}^{\star}$ 
be the word along $\mu$. Then $u_0= a^{b^n}v$ in $\QQ$. Killing $r$, $s$ and $t$,
retracts $\QQ$ onto the subpresentation
$\langle a, b \mid a^b = a^2 \rangle$.  
So $\bar{v} =  a^{b^n}$ in $\langle a, b \mid a^b = a^2 \rangle$, where $\bar{v} \in \set{ a^{\pm 1}}^{\star}$ is $v$ with all letters $t^{\pm 1}$ removed.  It follows that the exponent sum of the letters in $\bar{v}$ is $2^n$ and hence that $\mu$ has the desired property.  

Killing all generators other than $t$ defines a retraction $\phi$
of $\QQ$ onto $\langle t \rangle \cong \Z$ that is distance decreasing 
with respect to word metrics. But the image of $\phi \circ \pi: \Delta \to \Z$ has diameter at least
$2^n$ on account of $\mu$.  So $\IDiam(w_n) \geq 2^n$, as claimed.

\ms

It is easy to check that the van~Kampen diagram $\hat{\Delta}_n$
for $w_n$ constructed below admits a shelling down to their base
vertex that realises the bound $\FL(w_n) \preceq 2^n$.  So, as
$\IDiam(w_n) \leq \FL(w_n)$, we deduce that $\IDiam(w_n) \simeq
\FL(w_n) \simeq 2^n$.

The bound $\FFL(w_n) \preceq n$ follows from Proposition~\ref{crossing corridors}.  Nonetheless we will sketch a proof since the salient ideas, developing the cylinder example from Section~\ref{intro}, appear here more transparently than in the more general contexts of Section~\ref{aux}.  

The van~Kampen diagram $\hat{\Delta}_3$ for $w_3$ is depicted in Figure~\ref{Delta_n}.
The analogous construction of the diagram $\hat{\Delta}_n$ for $w_n$ should be clear.
Within $\hat{\Delta}_n$ there are four triangular
subdiagrams over the subpresentation $\langle a,t \mid [a,t]
\rangle$, in which strings of $a$-edges run vertically (in the sense of Figure~\ref{Delta_n}). 
Cut along each of these strings (except those of length $2^n$  at
the left and right sides of the diagrams) and insert back-to-back copies of the 
$\langle a,b \mid a^b=a^2 \rangle$-diagrams $\Omega_k$ that \emph{shortcut} $a^k$ to a word $u_k$ of length $\sim \log_2 k \preceq n$.  These \emph{shortcut diagrams} $\Omega_k$ are constructed in Proposition~\ref{BS shelling} and the way they are inserted is shown (in a more general context) in Figure~\ref{quarter fig}.   Call the resulting diagram $\Delta_n$.

\begin{figure}[ht]
\psfrag{a}{$a$}
\psfrag{b}{$b$}
\psfrag{T}{$s$}
\psfrag{t}{$t$}
\psfrag{tau}{$r$}
\psfrag{s}{$\star$}
\centerline{\epsfig{file=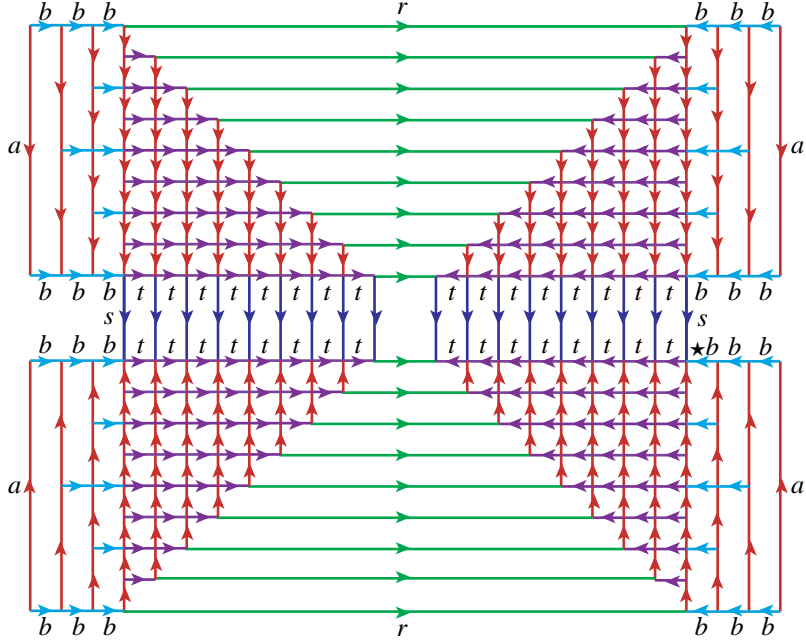}} \caption{The van~Kampen
diagram $\hat{\Delta}_3$ for $w_3$.} \label{Delta_n}
\end{figure}

For $0 \leq k \leq 2^n$ let $\hat{\rho}_k$ be the edge-paths in
$\hat{\Delta}_n$, forming concentric squares in
Figure~\ref{Delta_n}, labelled by $a^k s a^{-k} r^{-1} a^k
s^{-1} a^{-k} r$.  Next, for $0 \leq k \leq 2^n$ define  $\rho_k$ to be the edge-path in $\Delta_n$
that is obtained from $\hat{\rho}_k$ by replacing each subword $a^{\pm k}$ by its
shortcut ${u_k}^{\pm 1}$. (In particular $\rho_{2^n} = \partial \Delta_n$.) 
Note
that for all $k$, the length of $\rho_k$ is $\preceq n$.

We now briefly describe a full shelling of $\Delta_n$ that realises the
bound $\FFL(w_n) \preceq n$.  In the course of this shelling we encounter the
subdiagrams of $\Delta_n$ that have $\rho_k$ as their boundary loops.  The following lemma concerning the diagrams $\Omega_k$ of Proposition~\ref{BS shelling} is the key to shelling the subdiagram bounded by $\rho_{k+1}$ down to that bounded by $\rho_k$. During this shelling, the length of the boundary loop remains less than a linear function of $n$ because $\log_2 k \leq n$.  

\begin{lemma} \label{shelling lemma}
Let $\Pi$ be the van~Kampen diagram comprising a copy of $\Omega_{k+1}$ and a copy $\Omega_k$ joined to each side of a $t$-corridor along sides labelled $a^{k+1}$. Let $\star$ be a vertex on $\partial \Pi$ located at the start of either of the paths labelled $a^k$ along the sides of the $t$-corridor.  Then $\FL(\Pi, \star) \preceq \log_2 k$.
\end{lemma}

The proof of this lemma becomes clear when one considers  concurrently running the shellings of $\Omega_{k+1}$ and $\Omega_{k}$ of Proposition~\ref{BS shelling}, and a shelling of the $t$-corridor.

\ms 
We complete the proof of Theorem~\ref{Thm 1} by noting that the lower bound $\FFL(w_n) \succeq n$ is trivial as $\FFL(w_n) \geq \ell(w_n)=8n+8$, and so $\FFL(w_n) \simeq n$. 
\end{proof}

\bs 

\begin{proof}[Proof of Theorem~\ref{Thm 2}.]  
The lower bound of $2^n$ on $\IDiam(w_n)$ established above proves that $2^n  \preceq \IDiam(n)$.  So by (\ref{ineq1}) in Section~\ref{filling length} we see that $2^n \preceq \FL(n) \preceq \Area(n)$ and $\FFL(n) \preceq \Area(n)$ for all $n \in \N$. 

To establish $2^n \preceq \FFL(n)$ we show that the words $$w'_n \  :=  [s, \  a^{-  b^n} r a^{b^n}
r a^{b^n} r^{-1} a^{- b^n}]$$
satisfy $\FFL(w'_n) \geq 2^n$.  Figure~\ref{three corridors} shows a van~Kampen diagram for $w'_n$ built from a copy of the diagram $\hat{\Delta}_n$ for $w_n$, an $[r, s]$-2-cell, and a mirror image of $\hat{\Delta}_n$.  Thus $w'_n$ is null-homotopic.

\begin{figure}[ht]
\psfrag{a}{$a$}
\psfrag{b^n}{$b^n$}
\psfrag{a2n}{$a^{2^n}$}
\psfrag{t2n}{$t^{2^n}$}
\psfrag{at}{$(at)^{2^n}$}
\psfrag{T}{$s$}
\psfrag{t}{$t$}
\psfrag{tau}{$r$}
\psfrag{s}{$\star$}
\centerline{\epsfig{file=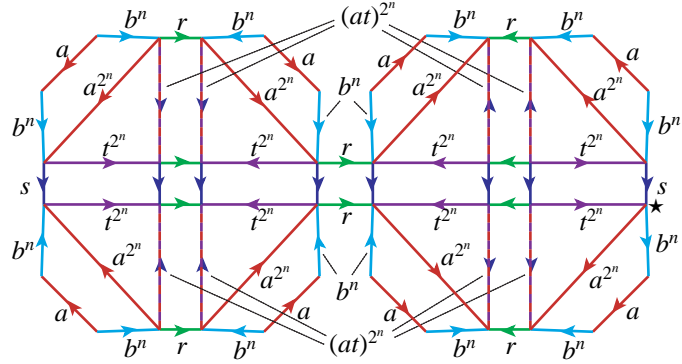}} \caption{A van~Kampen
diagram for $w'_n$.} \label{three corridors}
\end{figure}

To show that $\FFL(\Delta') \geq 2^n$ for \emph{all} van~Kampen diagrams $\Delta'$ for $w'_n$ we develop the argument used to establish the lower bound on diameter in the proof of Theorem~\ref{Thm 1}.  An $s$-corridor runs through $\Delta'$, and as there are three occurrences of $r$ in $w'_n$ and three of $r^{-1}$, each $r$ is joined to an $r^{-1}$ by an $r$-corridor that crosses the $s$-corridor.   This is illustrated in Figure~\ref{three more corridors}; note that the behaviour of the corridors could be more complex because each $r$-corridor could cross the $s$-corridor multiple times.  As shown in Figure~\ref{three more corridors}, let $u_1$, $u_2$ and $u_3$ be the words along the sides of the initial portions of these three $r$-corridors running from the first, second and third $r$ in $w'_n$ and ending where the corridor first meets the $s$-corridor.  (So  $u_1=u_3=(at)^{2^n}$ and $u_2$ is the empty word in the example of Figure~\ref{three corridors}.)

Retracting $\QQ$ onto the subpresentation
$\langle a, b \mid b^{-1} a b = a^2 \rangle$ by killing $r$, $s$ and $t$ we see that the exponent sum of the letter $a$, and hence also of $t$, in $u_1$ and in $u_3$ is $2^n$, while the exponent sum in $u_2$ is $0$.  The word along the side of the $s$-corridor is of the form $v_1r v_2 r v_3 r v_4$ where $v_i \in \set{r^{\pm 1}, t^{\pm 1}}^*$ and $v_i$ runs to the vertex at the end of $u_i$ for $i=1,2,3$ (see Figure~\ref{three more corridors}).   By considering the retraction of $\QQ$ onto $\langle t \rangle$ in which $a,b,r,s$ are killed, we see that the exponent sum of $t$ in $v_i$ is $2^n$ for $i=1,4$ and is $-2^n$ for $i=2,3$.  

Suppose $\mathcal{S}=(\Delta'_i)$ is a full shelling of $\Delta'$, in the course of which the base vertex is not required to be kept fixed.    
The retraction $\phi$ in which all generators other than $t$ are killed, defines a distance decreasing map $\phi$ of $\QQ$ onto $\langle t \rangle \cong \Z$.   And the edge-circuit $\phi \circ \pi \restricted{\partial \Delta'_i} : \partial \Delta'_i \to \Z$ has length at most $\ell(\partial \Delta'_i)$.
There are natural combinatorial maps $\psi_i: \Delta'_i \to \Delta'$ (only prevented from being injective by \emph{1-cell expansion} moves) under which $\psi_i(\partial  \Delta'_i)$ forms a contracting sequence of edge-circuits. Let $i$ be the least integer such that $\psi_i(\partial \Delta'_i)$ includes either the vertex at the end of $v_1$ or at the start of $v_4$ -- these are ringed by small circles in Figure~\ref{three more corridors}.  We will explain why $\ell(\partial \Delta'_i)$ is at least $2^n$ when this vertex $x$ is at the end of $v_1$.  A similar argument will show the same result to hold when $x$ is the vertex at the start of $v_4$.  Some vertex $y$ on $v_4$ must be included in $\psi_i(\partial \Delta'_i)$ because the contracting edge-circuit cannot have yet crossed the vertex at the start of $v_4$.  So $$\FFL(w'_n) \ \geq \  \ell(\partial \Delta'_i) \ \geq \ \abs{ \phi(x)- \phi(y)} \  \geq  \ 2^n,$$ as required.

\begin{figure}[ht]
\psfrag{u_1}{$u_1$}
\psfrag{u_2}{$u_2$}
\psfrag{u_3}{$u_3$}
\psfrag{v_1}{$v_1$}
\psfrag{v_2}{$v_2$}
\psfrag{v_3}{$v_3$}
\psfrag{v_4}{$v_4$}
\psfrag{a}{$a$}
\psfrag{b^n}{$b^n$}
\psfrag{a2n}{$a^{2^n}$}
\psfrag{t2n}{$t^{2^n}$}
\psfrag{at}{$(at)^{2^n}$}
\psfrag{T}{$s$}
\psfrag{t}{$t$}
\psfrag{tau}{$r$}
\psfrag{s}{$\star$}
\centerline{\epsfig{file=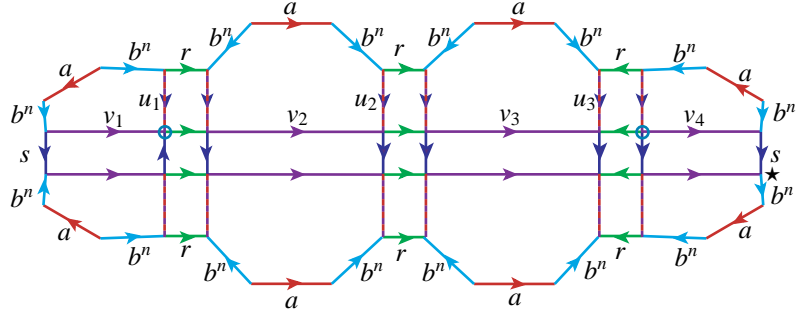}} \caption{A van~Kampen
diagram $\Delta'$ for $w'_n$.} \label{three more corridors}
\end{figure}

It is clear that $n \geq \FFFL(n)$ because for all null-homotopic words $w$ we have $\FFFL(w) \geq \ell(w)$.  Proposition~\ref{FFFL prop}, which is the culmination of a sequence of propositions in Section~\ref{aux}, will give $\FFFL(n) \simeq n$. Theorem~\ref{exp bound} will then give us $\Area(n) \preceq 2^n$ and the proof of the theorem will be complete.
\end{proof}

\section{Auxiliary propositions} \label{aux}

\ni In this section we provide a number  of results which build up to Proposition~\ref{FFFL prop} where we prove a linear upper bound on $\FFFL(n)$ in the presentation $\QQ$ of Theorem~\ref{Thm 1}. We begin in Proposition~\ref{BS shelling} with a technical result giving carefully controlled shellings of diagrams $\Omega_k$ that \emph{shortcut} $a^k$ in $\langle a,b \mid a^b=a^2 \rangle$ to a word of length $\preceq \log_2 k$.  Next Proposition~\ref{QQ_0 prop} claims the filling length of $\QQ_0  := \langle a,b,t \mid a^ba^{-2}, [a,t] \rangle$ admits a linear upper bound.  

Proposition~\ref{crossing corridors} gives a linear upper bound on $\FFL(w)$ for null-homotopic words in $\QQ$ that have exactly one pair of letters $r,r^{-1}$ and one pair $s,s^{-1}$, with the occurrences of the $r^{\pm 1}$ alternating with those of the $s^{\pm 1}$.  We show that such words have diagrams with one $s$- and one $r$-corridor with these corridors crossing only once.  These diagrams can be thought of as \emph{towers}, and are exemplified the diagrams $\Delta_n$ for $w_n$ constructed in the proof of Theorem~\ref{Thm 1}.    
In the proof of Proposition~\ref{crossing corridors} we refer back to Proposition~\ref{quarter} which provides controlled shellings for the four subdiagrams that \emph{tower} diagram breaks up into if we remove the $r$- and $s$-corridors.  

Next, in Proposition~\ref{2 Ts}, we establish an upper bound on $\FFFL(w)$ that is linear in $\ell(w)$, for null-homotopic words $w$ in $\QQ$ that have exactly one pair of letters $s,s^{-1}$.  The essential idea here is to find a diagram that can be fragmented into a number of subdiagrams, one for each $r$-corridor that crosses the $s$-corridor, and apply Proposition~\ref{crossing corridors} to each of these subdiagrams.  Proposition~\ref{no T} takes care of the case where there is no letter $s$ in $w$.  Finally, we prove Proposition~\ref{FFFL prop}: given a word that is null-homotopic  in $\QQ$ we construct a diagram that can be  fragmented into subdiagrams each of which contain at most one $s$-corridor, and then we apply Propositions~\ref{2 Ts} and \ref{no T}.

\begin{prop} \label{BS shelling}
Let $\PP$ be the presentation $\langle a,b \mid a^b=a^2
\rangle$.  Fix
$k>0$.  There is a word $u_k$ of length at most $12+4\log_2 k$ such
that $u_k=a^k$ in $\PP$.  Moreover, there is a $\PP$-van~Kampen diagram
$\Omega_k$ with boundary word $a^k{u_k}^{-1}$ \textup{(}as read from a
vertex $\star$\textup{)} satisfying the following condition:  if $\mu$ is the subarc of $\partial \Omega_k$
along which one reads $a^k$, then there is a shelling
$\mathcal{S}=\Omega_{k0}, \ldots, \Omega_{kp}$ that collapses $\Omega_k =
\Omega_{k0}$ to $\star= \Omega_{kp}$, such that each $\Omega_{ki}$ is a
subdiagram of $\Omega_k$ and, expressing the boundary circuit of $\Omega_{ki}$ as
$\mu_i$ followed by $\nu_i$, where $\mu_i$ is the maximal length
subarc of $\partial \Omega_{ki}$ that starts at $\star$ and follows
$\mu$, we have $\ell(\nu_i) \leq 12+4\log_2 k$.
\end{prop}

\begin{figure}[ht]
\psfrag{a}{$a$}%
\psfrag{b}{$b$}%
\psfrag{1}{$1$}%
\psfrag{2}{$2$}%
\psfrag{3}{$3$}%
\psfrag{4}{$4$}%
\psfrag{5}{$5$}%
\psfrag{6}{$6$}%
\psfrag{7}{$7$}%
\psfrag{8}{$8$}%
\psfrag{9}{$9$}%
\psfrag{10}{$10$}%
\psfrag{11}{$11$}%
\psfrag{12}{$12$}%
\psfrag{13}{$13$}%
\psfrag{14}{$14$}%
\psfrag{15}{$15$}%
\psfrag{star}{$\star$}%
\centerline{\epsfig{file=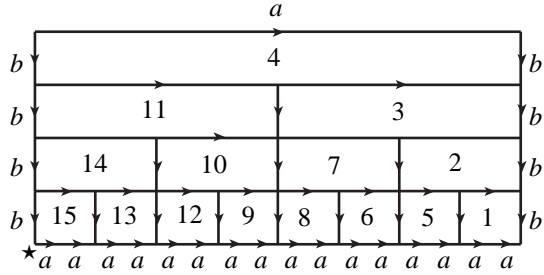}} \caption{The van~Kampen
diagram $\Xi$ for $a^{2^4}a^{- b^4}$.} \label{BS fig}
\end{figure}

\begin{proof}  Let $m$ be the least integer such that $2^m \geq k$.  So $m < 1 + \log_2 k$.  Let $\Xi$ be the standard van~Kampen diagram that demonstrates the equality $a^{b^m}=a^{2^m}$ and is depicted in Figure~\ref{BS fig} in the case $m=4$.  Let $(\Xi_i)$ be the shelling of $\Xi$ in which each $\Xi_i$ is a subdiagram of $\Xi$ and $\Xi_{i+1}$ is obtained from $\Xi_i$ as follows.  A 1-cell collapse is performed on a \emph{spike} of $\Xi_i$ if possible.  Otherwise, of the rightmost 2-cells in $\Xi_i$, let $e^2$ be the lowest in the sense of Figure~\ref{BS fig}.  Let $e^1$ be the right-most edge of the lower horizontal side of $e^2$.  Then $e_1 \in \partial \Xi_i$.  Perform a 2-cell collapse on $(e^2,e^1)$.  As an illustration, the numbers in the 2-cells in Figure~\ref{BS fig} show the order in which the 2-cells are collapsed.  

When there is no spike in $\Xi_i$, its anticlockwise boundary circuit, starting from $\star$, follows \emph{horizontal} (in the sense of Figure~\ref{BS fig}) edges labelled by $a$, then travels \emph{upwards} through the boundary of $m$ 2-cells (visiting at most three 1-cells of each) and then \emph{descends} back to $\star$ along edges labelled $b$.  Removing the initial horizontal path, the number of edges traversed is at most $4m$.  And, as the total length of the 1-dimensional portions of $\partial \Xi_i$ is at most $8$, and $m < 1 + \log_2 k$, we deduce that the length of $\partial \Phi_i$ minus the length of the initial horizontal arc, is at most $12 + 4 \log_2 k$.  

In the case $k=2^m$ we find that defining $\Omega_i := \Xi_i$ for all $i \geq 0$ gives the asserted result.  For the case $k \neq 2^m$, let $c$ be the maximum $i$ such that $a^k$ is a prefix of the word $w$ one reads anticlockwise around $\partial \Xi_i$, starting from $\star$.  Then defining $u_k:= {w_0}^{-1}$, where $w=a^kw_0$ and $\Omega_{ki} := \Xi_{c+i}$ for all $i \geq 0$, we have our result.   \end{proof}

\ms

\ms

\begin{prop} \label{QQ_0 prop}
The filling length function $\FL:\N \to \N$ of 
\begin{eqnarray*}
\QQ_0 \ & := & \ \langle a,b,t \mid a^ba^{-2}, [a,t]
\rangle
\end{eqnarray*}
admits a linear upper bound.
\end{prop}

\begin{proof}
Corollary~E1 of \cite{BGSS} says that an HNN extension of a finitely generated free group with finitely many stable letters, in which the associated subgroups are all finitely generated, is asynchronously automatic. This applies to $\QQ_0$.     Theorem~3.1 in \cite{Gersten5} says (in different language) that if a group is asynchronously combable then its filling length function admits a linear upper bound.
\end{proof}

\begin{figure}[ht]
\psfrag{a^k}{$a^k$}%
\psfrag{u_k}{$u_k$}%
\psfrag{t^{j-k}}{$t^{j-k}$}%
\psfrag{at}{$(at)^k$}%
\psfrag{t^j}{$t^j$}%
\psfrag{t^k}{$t^k$}%
\psfrag{t^j-k}{$t^{j-k}$}%
\psfrag{w}{$w$}%
\psfrag{D_i}{$\Delta_i$}%
\psfrag{Dp}{$\Delta'$}%
\psfrag{Dpp}{$\Delta''$}%
\psfrag{D_m}{$\Delta_m=\Lambda$}%
\psfrag{D_0}{$\Delta=\Delta_0$}%
\psfrag{leqL}{$\leq L$}%
\centerline{\epsfig{file=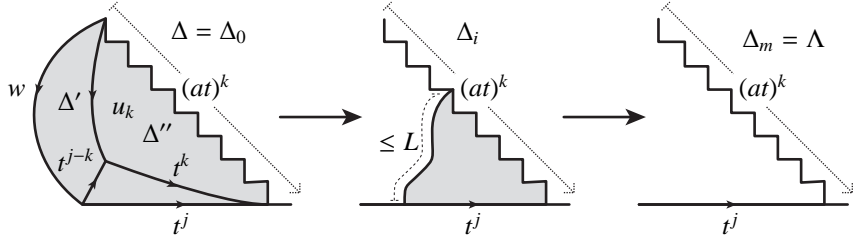}} \caption{Shelling $\Delta$ down to $\Lambda$.} \label{quarter_shelling}
\end{figure}

\begin{prop} \label{quarter}
Suppose $(at)^{-k}wt^j$ is a null-homotopic word in $\QQ_0= \langle a,b,t \mid a^ba^{-2}, [t,a] \rangle.$  Let $\Lambda$ be the 1-dimensional van~Kampen diagram for $$t^j(at)^{-k}(at)^{k-1} at^{-(j-1)}$$ constructed by assembling 1-cells in $\mathbb{R}^2$ as depicted at the right in Figure~\ref{quarter_shelling}.  There is a $\QQ_0$-van~Kampen diagram $\Delta$ for  $(at)^{-k}wt^j$ with the following properties. There is a shelling of $\Delta$ through a sequence of diagrams $\Delta=\Delta_0, \Delta_1, \ldots, \Delta_m=\Lambda$ with the portion $t^j (at)^{-k}$ of $\partial \Delta$ left undisturbed throughout \textup{(}see Figure~\ref{quarter_shelling}\textup{)}.  Let $\nu_i$ be the maximal length arc of the boundary circuit of $\Delta_i$ contained entirely in $\Lambda$.  There exists $C>0$, depending only on $\QQ_0$ such that $L:= \max_i ( \ell(\partial \Delta_i)- \ell(\nu_i)) \leq C \ell(w)$.  
\end{prop}

\begin{figure}[ht]
\psfrag{a^k}{$a^k$}
\psfrag{a^{k-1}}{$a^{k-1}$}
\psfrag{a^{k-2}}{$a^{k-2}$}
\psfrag{a^{k-3}}{$a^{k-3}$}
\psfrag{O_6}{$\Omega_6$}
\psfrag{O_4}{$\Omega_4$}
\psfrag{O_3}{$\Omega_3$}
\psfrag{O_2}{$\Omega_2$}
\psfrag{O_1}{$\Omega_1$}
\psfrag{u_6}{$u_6$}
\psfrag{u_4}{$u_4$}
\psfrag{u_3}{$u_3$}
\psfrag{u_2}{$u_2$}
\psfrag{u_1}{$u_1$}
\psfrag{u_k}{$u_k$}
\psfrag{u_k-1}{$u_{k-1}$}
\psfrag{u_k-2}{$u_{k-2}$}
\psfrag{u_k-3}{$u_{k-3}$}
\psfrag{t^{j-k}}{$t^{j-k}$}
\psfrag{t^j}{$t^j$}
\psfrag{t^k}{$t^k$}
\psfrag{w}{$w$}
\psfrag{a}{$a$}
\psfrag{t}{$t$}
\centerline{\epsfig{file=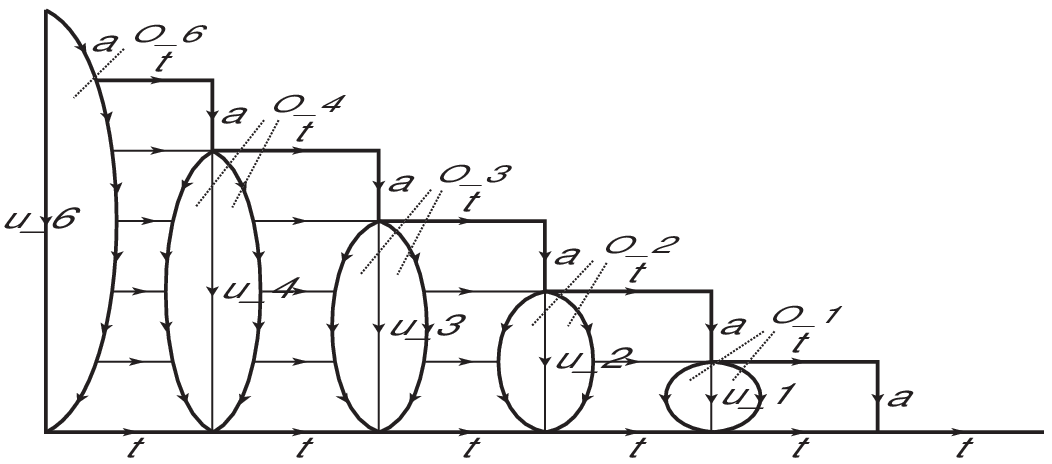}} \caption{$\Delta''$.} \label{quarter fig}
\end{figure}

\begin{proof}
We consider first the case $k \geq 0$.
In $\QQ_0$ we find that $w^{-1}=t^j (at)^{-k} = t^{j-k} a^{-k} = t^{j-k}{u_k}^{-1}$, where $u_k$ is the word  of Proposition~\ref{BS shelling} that has length at most $12+4 \log_2 k$.  These equalities are displayed in the left-most diagram of Figure~\ref{quarter_shelling}, which shows the framework of the van~Kampen $\Delta$: a union of a diagram $\Delta'$ for $t^{j-k}{u_k}^{-1} w$, a diagram $\Delta''$ for $(at)^{-k} u_k t^k$ and a tripod; the lower triangular region in the figure folds up to give a tripod, the exact configuration of which depends on the relative signs of $j$, $k$ and $j-k$.  

Note that $j-k$ is the exponent sum of the $t^{\pm 1}$ in $w$ as $\QQ_0$ retracts onto $\langle t \rangle \cong \Z$, and so $\abs{j-k} \leq \ell(w)$.  And $k\leq 2^{\ell(w)}$ because killing $t$ retracts $\QQ_0$ onto $\PP$.  It follows that $\ell(u_k)\leq 12+4 \ell(w)$ and  $\ell(t^{j-k}{u_k}^{-1} w) \leq 2 \ell(w) + \ell(u_k) \leq 12+6 \ell(w)$.  By Proposition~\ref{QQ_0 prop} we can take $\Delta'$ to be a van~Kampen diagram for $t^{j-k}{u_k}^{-1} w$ with filling length at most a constant times $12+6 \ell(w)$.  We can cut along the edge-path in $\Delta$ labelled by $u_kt^{-(j-k)}$, leaving $\Delta'$ attached to the rest of the diagram at only one vertex, and then shell $\Delta'$, and in the process the length of the non-$t^j(at)^{-k}$-portion of the boundary curve has length at most a constant times $\ell(w)$.  

The word $t^{k}(at)^{-k} a^k$ admits an obvious diagram with vertical $t$-corridors (as shown in Figure~\ref{quarter fig}) of height $k-1, k-2, \ldots, 1$.  We cut along the vertical paths labelled by powers of $a$ and insert copies of the diagrams $\Omega_i$ of Proposition~\ref{BS shelling} and their mirror images, as shown in Figure~\ref{quarter fig} (illustrated in the case $k=6$).  The shellings of Lemma~\ref{shelling lemma} can be composed to give a shelling down to $\Lambda$ that realises the asserted bound on $L$.   \end{proof}

For a word $w$ we use $\ell_t(w)$ to denote the number of letters $t^{\pm1}$ in $w$.   

\begin{prop} \label{crossing corridors}
Suppose $w$ is a null-homotopic word over 
$$\QQ:=\langle \ a,b, r, s, t \ \mid \ a^ba^{-2}, [t,a], [r, at], [r,s], [s,t]  \ \rangle$$
and $\ell_{r}(w)=\ell_{s}(w)=2$ with the occurrences of $r^{\pm 1}$ and $s^{\pm 1}$ alternating in $w$.  Then $\FFL(w) \leq C \ell(w)$ where $C>0$ depends only on $\QQ$.
\end{prop}

\begin{proof}
In any van~Kampen diagram $\Delta$ for $w$ there is one $r$-corridor and one $s$-corridor.  The condition that the occurrences of $r^{\pm 1}$ and $s^{\pm 1}$ in $w$ alternate is equivalent to saying that these corridors must cross at least once in $\Delta$.  
 
Let $\Delta$ be a minimal area diagram for $w$.  Suppose, for a contradiction, that the $s$-corridor and $r$-corridor in $\Delta$ cross more than once.  Then there is a subdiagram between the two corridors with boundary word $w_0 =uv$ where $u \in \set{t^{\pm 1}}^{\star}$ and $v \in \set{(at)^{\pm 1}}^{\star}$.  Both $u$ and $v$ must be reduced words as otherwise $\Delta$ would not be a reduced diagram and hence not be of minimal area.  Killing $r$, $s$ and $t$ retracts $\QQ$ onto a group in which $a$ has infinite order, and so $v$ is the empty word.  Killing all the generators other than $t$ retracts $\QQ$ onto $\langle t \rangle =\Z$ and so $u$ is the empty word.   It follows that there are cancelling $[r,s]$-2-cells where the two corridors cross and so $\Delta$ is not of minimal area -- a contradiction.

A similar argument shows that  $\Delta$  contains no $r$- or $s$-annuli.    

Conclude that $\Delta$ consists of an $r$-corridor, an $s$-corridor and four subdiagrams of the form where Proposition~\ref{quarter} applies.  Produce a new van~Kampen diagram $\Delta'$ for $w$ by replacing the four subdiagrams that minimise the length $L$ of Proposition~\ref{quarter}.  A shelling of $\Delta'$ realising the asserted bound is obtained by running shellings of the four subdiagrams and the two corridors concurrently in the obvious way so that the diagram is eventually shelled to the $[r,s]$-2-cell, and then to a single vertex.   \end{proof}

\begin{prop} \label{2 Ts}
Suppose $w$ is a null-homotopic word in $\QQ$ \textup{(}defined above\textup{)} and $\ell_s(w)=2$.  Then $\FFFL(w)\leq C\ell(w)$ where $C>0$ depends only on $\QQ$.  
\end{prop}

\begin{proof}
Let $\Delta$ be a reduced van~Kampen diagram for $w$. 
We will use the layout of the $r$-corridors and the one $s$-corridor in $\Delta$ as a template for the construction of another van~Kampen diagram $\Delta_2$ for $w$ that will admit a shelling realising the asserted bound. 

Suppose $\mathcal{C}$ is an $r$-corridor in $\Delta$ that does not cross the $s$-corridor. 
The word $w_{\mathcal{C}}$ along the sides of $\mathcal{C}$ is in $\set{ (at)}^{\star}$ and is reduced because $\Delta$ is a reduced diagram,  and so must be $(at)^k$ for some $k \in \Z$.    Killing all defining generators other than $t$ retracts $\QQ$ onto $\langle t \rangle$. So $k \leq \ell(\lambda)$, where $\lambda$ is a portion of the boundary circuit of $\Delta$ connecting the end points of a side of $\mathcal{C}$.  

If follows that if we remove any number of $r$-corridors that do not cross the $s$-corridor from  $\Delta$, then the length of the boundary circuit of each connected component is at most $2\ell(w)$.  

Suppose we remove \emph{all} of the $r$-corridors that do not cross the $s$-corridor from $\Delta$.  Define $\Delta_0$ to be the connected component that contains the $s$-corridor.   
All of the other connected components have boundary words that are null-homotopic in $$\QQ_0= \langle a,b,t \mid a^ba^{-2}, [t,a] \rangle,$$ which is both a retract and a subpresentation of $\QQ$.  
Obtain $\Delta_1$ from $\Delta$ by replacing all these subdiagrams by $\QQ_0$-diagrams of minimal $\FFL$.  

Repeating the following gives a shelling of $\Delta_1$ down to $\Delta_0$ in the course of which the boundary circuit has length at most $C_0 \ell(w)$, where $C_0$ is a constant that depends only on $\QQ$.  Choose an $r$-corridor $\mathcal{C}$ in $\Delta_1$ such that, of the two components we get by removing $\mathcal{C}$, that which does not contain the $s$-corridor contains no $r$-corridor.  Cut along one side of $\mathcal{C}$ using 
\emph{1-cell expansion} moves, and one \emph{fragmentation} move.  Next use \emph{1-cell collapse} and \emph{2-cell collapse} moves to remove the 2-cells along $\mathcal{C}$. By the remarks above, both connected components have boundary circuits of  length at most $2 \ell(w)$.  Collapse the component that does not contain the $s$-corridor down to a single vertex using a minimal $\FFL$ shelling, in the course of which the boundary circuit has length at most a constant times $2\ell(w)$ by Proposition~\ref{QQ_0 prop}.  (In fact, a shelling down to $\Delta_0$ within the required bound, can be achieved without the fragmentation move if care is taken over basepoints.)   

It remains to show that the word $w_0$ around $\partial \Delta_0$ admits a van~Kampen diagram with a full, fragmenting, free shelling in which the sum of the lengths of the boundaries of the components are at most a constant times $\ell(w_0)$.  For then we can take $\Delta_2$ to be $\Delta_1$ with $\Delta_0$ replaced by this diagram.  

First suppose  $\ell_{r}(w_0) = 0$.  Then $w_0$ is null-homotopic in the retract $\QQ_1:=\langle \ a,b, s, t \ \mid \ a^ba^{-2}, [t,a], [s,t] \ \rangle$, and the length of the $s$-corridor in any reduced $\QQ_1$-diagram $\Delta'_0$ for $w_0$ is at most $\ell(w_0)/2$ on account of the retraction onto $\langle t \rangle$.  Assume that the two components of $\Delta'_0$ we get on removing the $s$-corridor are $\QQ_0$-diagrams of minimal $\FFL$.  Then we can collapse $\Delta'_0$ by shelling each of these components and the $s$-corridor in turn, and using Proposition~\ref{QQ_0 prop} it is easy to check that the length of the boundary circuit remains at most a constant times $\ell(w_0)$.

Next suppose $\ell_{r}(w_0) = 2$. Then Proposition~\ref{crossing corridors} applies and gives us the result we need.  

\begin{figure}[ht]
\psfrag{tau}{$r$}
\psfrag{taue}{$r^{\varepsilon}$}
\psfrag{tau-e}{$r^{-\varepsilon}$}
\psfrag{w_1}{$w_1$}
\psfrag{C}{$\mathcal{C}$}
\psfrag{T}{$s$}
\centerline{\epsfig{file=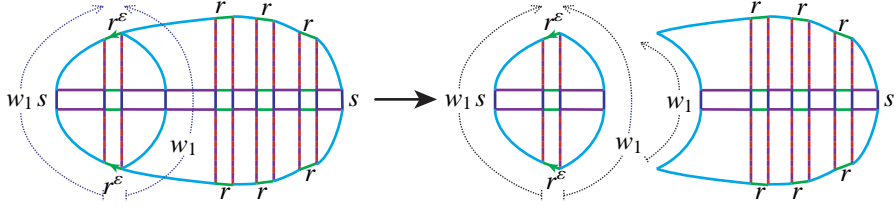}} \caption{Shelling away one $r$-corridor.} \label{frag cor}
\end{figure}

Finally, suppose $\ell_{r}(w_0) > 2$.  Then there is a subword $r^{\varepsilon} w_1 r^{-\varepsilon}$ in $w_0$, where $\varepsilon = \pm 1$, $\ell_{r}(w_1)=0$ and $\ell_{s}(w_1)=1$.  As $r^{\varepsilon} w_1 r^{-\varepsilon}= w_1$ in $\QQ$, there is a $\QQ$-van~Kampen diagram for $w_0$ that we can shell by cutting along an edge-path labelled by $w_1$ to cut the diagram into two, as shown in Figure~\ref{frag cor}, and the shelling the two components.  One of these components is a diagram for $r^{\varepsilon} w_1 r^{-\varepsilon} {w_1}^{-1}$, and this we shell first as per Proposition~\ref{crossing corridors}.  The remaining component has boundary length $\ell(w_0)-2$ and  includes two fewer letters $r^{\pm 1}$, and so by continuing inductively we can find a shelling for 
which $\FFFL$ is at most a constant times $\ell(w_0)$. 
\end{proof}

\begin{prop} \label{no T}
Suppose $w$ is a null-homotopic word in $$\QQ_2:=\langle \ a,b, r, t  \mid  a^ba^{-2}, [t,a], [r, at]  \ \rangle.$$  Then $\FL(w)\leq C\ell(w)$ where $C>0$ depends only on $\QQ_2$.  
\end{prop}

\begin{proof}
The method used in the proof of Proposition~\ref{crossing corridors}, to reduce to the case where all the $r$-corridors cross the $s$-corridor, gives this result. 
\end{proof}

\begin{prop} \label{FFFL prop}
Suppose $w$ is a null-homotopic word in $\QQ$.  Then $\FFFL(w)\leq C\ell(w)$ where $C>0$ depends only on $\QQ$.  
\end{prop}

\begin{proof}
The cases where $\ell_s(w)=0$ and $\ell_s(w)=2$ are dealt with by Propositions~\ref{no T} and \ref{crossing corridors}, respectively.   For the case $\ell_s(w)>2$ we take a similar approach to that used in the proof of Proposition~\ref{crossing corridors} to control $\FFFL(w_0)$.  

There is a subword $s^{\varepsilon} w_1 s^{-\varepsilon}$ in $w$, where $\varepsilon = \pm 1$ and $s^{\varepsilon} w_1 s^{-\varepsilon}=w_1$ in $\QQ$.  So we can find a $\QQ$-van Kampen diagram for $w$ that can be severed into two components, one of which is a diagram for $s^{\varepsilon} w_1 s^{-\varepsilon}{w_1}^{-1}$, and the other of which is a diagram for a word of length $\ell(w)-2$ that has two fewer letters $s^{\pm 1}$.  The former of these two components can be shelled as per Proposition~\ref{crossing corridors}.  Continuing inductively we see that the other component can be taken to be a diagram that admits a shelling in which the boundary circuit has length at most a constant times $\ell(w)$.  
\end{proof}

\section{Quasi-isometry invariance} \label{qi}

\ni In this section we prove Theorem~\ref{qi thm}, the quasi-isometry invariance amongst finitely presented groups of the functions $\FL$, $\FFL$ and $\FFFL$ up to $\simeq$-equivalence.  Our approach is to monitor how filling length, in its three guises, behaves in the standard proof that finite presentability is a quasi-isometry invariant \cite[page 143]{BrH}; other quantified versions of this proof exist in the literature (the first being \cite{Alonso}, where the focus is on $\Area$), so our exposition here will be brief.  

We have quasi-isometric groups $\Gamma$ and $\Gamma'$ with finite presentations $\PP= \langle \AA \mid \RR \rangle$ and $\PP'= \langle \AA' \mid \RR' \rangle$, respectively.  So there is a quasi-isometry $f: (\Gamma,d_{\AA}) \to (\Gamma',d_{\AA'})$ 
with quasi-inverse $g: (\Gamma',d_{\AA'}) \to (\Gamma,d_{\AA})$.
We begin by showing $\FL_{\PP} \simeq \FL_{\PP'}$. 

Suppose $\rho'$ is an edge-circuit in the Cayley graph
of $\PP'$, visiting vertices $v_0, v_1, \dots, v_n$
in order with $v_0=v_n$. Consider a circuit $\rho$ in the Cayley graph $\Cay^1(\PP)$ 
obtained by choosing geodesics joining $g(v_i)$ to $g(v_{i+1})$ for $i=0,\ldots, n-1$.  Note that $\ell(\rho)$ is at most a constant times  $\ell(\rho')$. 
Fill $\rho$ with a van~Kampen diagram $\Delta$ over $\PP$ admitting a shelling $\mathcal{S}$ of filling length $\FL_{\PP}(\ell(\rho))$ and  
use $f$ to map $\Delta^{(0)}$ to $\Gamma'$. Then join $f(a)$ to $f(b)$ by a geodesic whenever $a$ and $b$ are the end points of an edge in $\Delta$.  The result is a combinatorial map $\pi'_0 : {\Delta'_0}^{(1)} \to \Cay^1(\PP')$ filling a loop $\rho'_0$, where $\Delta'_0$ is obtained from $\Delta$ by subdividing each of its edges into edge-paths of length at most some constant.  

Interpolate between  $\rho'$ and $\rho'_0$ by joining $v_i$ to $f(g(v_i))$ for every $i$, to build a map $\pi'_1 : {\Delta'_1}^{(1)} \to \Cay^1(\PP')$ where $\Delta'_1$ is obtained from $\Delta'_0$ by attaching an annulus $A$ of $n$ 2-cells around the boundary.   
     
One obtains
a shelling $\mathcal{S}'_1$ of $\Delta'_1$ down to the base vertex with $\ell(\mathcal{S}'_1)$ at most a constant times $(1+\FL_{\PP}(\ell(\rho)))$ by first shelling away $A$ to leave just a stalk from $v_0$ to $f(g(v_0))$, and then running a shelling  of $\Delta'_0$ modelled on $\mathcal{S}$; whenever $\mathcal{S}$ demands the collapse of a 1-cell in $\Delta$, one collapses all the 1-cells in the corresponding edge-path in $\Delta'_0$.     

Although the words labelling the 2-cells of $\Delta'_1$  are null-homotopic they may fail to be in $\RR'$, so $\Delta'_1$ is not yet  a van~Kampen diagram over $\PP'$.  To rectify this one should replace the 2-cells of $\Delta'_1$ by van~Kampen diagrams over $\PP'$, each of area at most some uniform constant. But a problem arises in that van~Kampen diagrams can be singular 2-discs, so gluing them in place of 2-cells may destroy planarity.  
One gets around this by replacing the 2-cells of $\Delta'_1$ one at a time
in the following manner. If the boundary circuit
of the 2-cell $e^2$ that is to be replaced is not embedded, then 
we focus on the innermost embedded circuit $\sigma$ in the 1-skeleton that
encloses a disc containing $e^2$ (this has length less than the boundary
circuit of $e^2$). We delete the entire subdiagram enclosed by $\sigma$ and replace it with a van Kampen for the word labelling $\sigma$.
The result is a van Kampen diagram $\Delta'_2$ for $\rho'$ over $\PP'$.  We obtain a shelling $\mathcal{S}'_2$ for $\Delta'_2$ by altering $\mathcal{S}'_1$: each time we discarded some connected component of the set of edges inside some $\sigma$ we contract it (more strictly, its pre-image) and all the 2-cells it encloses to a single vertex in every one of the diagrams comprising the shelling, and each time we fill some 2-cell $e^2$ with a van~Kampen diagram $D$ we shell out all of $D$ when we had been due to perform a 2-cell collapse move on $e^2$.  [This could fail to give a shelling when a 2-cell collapse move in $\mathcal{S}'_1$ removes a pair $(e^2,e^1)$, where $e^1$ is one of the now contracted edges -- but then a 1-cell expansion followed by a 2-cell collapse producing the same effect can be used instead.]
The difference between $\ell(\mathcal{S}'_1)$ and $\ell(\mathcal{S}'_2)$ is then at most some additive constant.  Deduce that $\FL_{\PP'} \preceq \FL_{\PP}$.  Interchanging $\PP$ and $\PP'$ we have $\FL_{\PP} \preceq \FL_{\PP'}$ and so $\FL_{\PP} \simeq \FL_{\PP'}$.

The proof that  $\FFL_{\PP} \simeq \FFL_{\PP'}$ is essentially the same, except we consider free shellings, we discard the stalk between  $v_0$ and $f(g(v_0))$, and we replace $\FL_{\PP}(\ell(\rho))$ by $\FFL_{\PP}(\ell(\rho))$.  To show $\FFFL_{\PP} \simeq \FFFL_{\PP'}$ we additionally allow free and fragmenting shellings and we use $\FFFL_{\PP}(\ell(\rho))$ in place of $\FL_{\PP}(\ell(\rho))$;  no further technical concerns arise.

\section{Riemannian versus combinatorial filling length} \label{rc}

\ni Suppose $c: [0,1] \to X$ is a loop in a metric space $X$.   

A \emph{based null-homotopy} $H$ of $c$ is a continuous map $H : [0,1]^2 \to X$ for which $H(0,t)=H(1,t)=c(0)$ for all $t$ and, defining $H_t : [0,1] \to X$ by $H_t(s) = H(s,t)$, we have $H_0 = c$ and $H_1(s) = c(0)$ for all $s$. 

A \emph{free null-homotopy} $H$ of  $c$  is a continuous map $H : [0,1]^2 \to X$ such that $H(0,t)=H(1,t)$ for all $t$, and   $H_0$ and $H_1$ are $c$ and a constant function, respectively.

Let $\mathcal{S}$ be the set of subspaces $S$ of $[0,1]^2$ that have $[0,1] \times \set{0} \subset S$ and are the union of a finite family of  closed  triangular discs $\Delta_i = [(a_i,0), (b_i,0), (c_i,d_i)]$ with $c_i \in [a_i, b_i]$ and $d_i \in (0, 1]$.  The fibres $S_t$ of the projection mapping points in $S$ to their second co-ordinate are empty or are disjoint unions $\bigsqcup_{i=1}^{k_t} I_{t,i} \times \set{t}$ of closed intervals $I_{t,i}$. At finitely many critical $t$-values $\tau_l$, some of the intervals comprising the fibre bifurcate or collapse and vanish.

A \emph{free and fragmenting null-homotopy} $H$ of  $c$  is a continuous map $H :  S \to X$ for some $S \in \mathcal{S}$ where, defining $H_t$ to be the restriction of $H$ to $S_t$, we find $H_0 = c$ and $H_t(x) = H_t(y)$ for all $t$ whenever $x$ and $y$ are the end points of some $I_{t,i}$.  We define $\ell(H_t)$ to be the sum of the lengths of the $k_t$ loops in $X$ defined by $H_t$.    Taking $S$ to be a single triangle reduces to the case of a free null-homotopy.   

In each of the three settings above define $\ell(H)  \ := \ \sup_{t \in [0,1]} \ell(H_t)$,
and then
\begin{eqnarray*}
\FL(c) & := &  \inf \set{  \ \ell(H) \mid \textup{ based null-homotopies } H \textup{ of } c  \ } \\ 
\FFL(c) & := &  \inf \set{  \ \ell(H) \mid \textup{ free null-homotopies } H \textup{ of } c \ } \\ 
\FFFL(c) & := &  \inf \set{ \ \ell(H) \mid \textup{ free and fragmenting null-homotopies } H \textup{ of } c \ }. 
\end{eqnarray*}
For $\textup{M}= \FL, \FFL$ or $\FFFL$, define $\textup{M}_X: [0, \infty) \to [0, \infty)$ by
$$\textup{M}_X(l) \ = \ \sup \left\{ \left. \  \textup{M}(c) \ \right| \textup{ null-homotopic loops } c \textup{ with } \ell(c) \leq l \ \right\}.$$ 

The following lemma gives sufficient conditions for $\FL_X$, $\FFL_X$ and $\FFFL_X$ to be well-defined --- conditions enjoyed by the universal cover of any closed connected Riemannian manifold, for example.  

\begin{lemma}   \label{well defined}  
Suppose $X$ is the universal cover of a compact geodesic space $Y$ for which there exist $\mu, L > 0$ such that every loop of length less than $\mu$ admits a based null-homotopy of filling length less than $L$.   Then $\FL_X$, $\FFL _X$ and $\FFFL _X$ are well-defined functions $[0,\infty)\to [0,\infty)$.
\end{lemma}

\begin{proof} 
The proof of Lemma~2.2 in \cite{BR1} can readily be adapted to this context.  
In brief, we first show that every rectifiable loop $c$ in $X$ admits a based null-homotopy with finite filling length -- apply a compactness argument to an arbitrary based null-homotopy for $c$ to partition $c$ into finitely many loops of length at most $\mu$; by hypothesis each such loop has finite filling length and it follows that $c$ has finite filling length.  

Next  suppose $c$ has length $l$ and assume (by reducing $\mu$ if necessary) that balls of radius $\mu$ in $Y$ lift to $X$.
Cover $Y$ with a maximal collection of disjoint closed balls of radius $\mu/5>0$;  let $\Lambda \subset X$ be the set of lifts of their centres.  Subdivide $c$ into $m \leq 1+10\ell/\mu$ arcs with end-points $v_i$, each of length at most $\mu/10$; each $v_i$ lies within $\mu/5$ of some $u_i \in \Lambda$; form a piecewise geodesic loop $c'$ approximating $c$ by connecting-up these $u_i$.  Loops made up of the portion of $c$ from $v_i$ to $v_{i+1}$ and geodesics $[u_i,u_{i+1}]$, $[u_i, v_i]$ and   $[u_{i+1}, v_{i+1}]$ have length at most $\mu$.  Homotopy discs for these loops together form a \emph{collar} between $c$ and $c'$.
By passing across these discs one at a time, it is possible to homotop $c$ across the collar to a loop made up of $c'$ and a stalk of length $\mu / 5$ from $c(0)$ to $u_0$, encountering loops only of length at most a constant (depending on $L$ and $\mu$)  times $l$ en route.       
Modulo the action of $\pi_1Y$,  there are only finitely many such piecewise geodesic loops such as $c'$ and, by our earlier argument, each one admits a filling of finite filling length.  It follows that  
$\FL_X$, and hence $\FFL_X$ and $\FFFL_X$, are well-defined functions. 
\end{proof}

\ms 
\begin{proof}[Proof of Theorem~\ref{translate}.]  
Fix a basepoint $p \in X$.  Define a quasi-isometry $\Phi$ mapping the Cayley graph of $\PP = \langle \AA \mid \RR \rangle$ to $X$ by  choosing a geodesic from $p$ to its translate $a\cdot p$ for each  $a \in \AA$, and then extending equivariantly.  Let $\Psi$ be a quasi-isometry from $X$ to $\Gamma$ sending $x \in X$ to some $\gamma$ such that $\gamma.p$ is a point of $\Gamma.p$ closest to $x$.  

A path in $X$ is called {\em word-like} (following \cite{Bridson6}) if it is the image in $X$ of an edge-path in the Cayley graph.  For each $r \in \RR$, let $c_r$ denote the word-like loop in $X$, based at $p$ that is the image of an edge-circuit in the Cayley graph labelled $r$. Map the Cayley 2-complex of $\PP$ to $X$ by choosing a disc-filling arising from a based null-homotopy of   
finite filling length for each $c_r$.   

We will show first that $\FL_{X} \preceq \FL_{\PP}$, $\FFL_{X} \preceq \FFL_{\PP}$ and $\FFFL_{X} \preceq \FFFL_{\PP}$.  As in the proof of Lemma~\ref{well defined}, a collar between an arbitrary rectifiable loop $c$ in $X$ and a word-like loop $c'$, can be used to show there is no change in the $\simeq$ classes of $\FL_{X}$, $\FFL_{X}$ or $\FFFL_{X}$ if one takes the suprema in their definitions to be over fillings of word-like loops only: for $\FL_{X}$ one notes that $c$  can be homotoped across the collar  to a loop based at $c(0)$ that is obtained from $c'$ by attaching a stalk from $c(0)$ to $c'(0)$, and one need pass through loops of length no more than $C\ell(c)+C$ en route, where $C$ is a constant independent of $c$; for  $\FFL_{X}$ and  $\FFFL_{X}$, the stalk is abandoned and the homotopy is between $c$ and $c'$.

One gets an upper bound on the filling length of a word-like loop $c$ in $X$  by taking the image in $X$ of a minimal filling length van~Kampen diagram $\Delta$.  The progress of the boundary circuit in the course of a shelling of $\Delta$ dictates a sequence of stages in a null-homotopy of $c$.  Using Lemma~\ref{well defined}, we can interpolate between these stages in a way that increases the length of the curve by no more than an additive constant, and so we get $\FL_{X} \preceq \FL_{\PP}$.  The proof that $\FFL_{X} \preceq \FFL_{\PP}$ and $\FFFL_{X} \preceq \FFFL_{\PP}$ can be completed likewise.

Now we address $\FL_{\PP} \preceq \FL_X$.
Consider a word-like loop $c: [0,1] \to X$ corresponding to a null-homotopic word $w$ over $\PP$ of length $n$.  Fix a constant $\lambda > \max_{a \in \AA} d_X(p, a.p)$. Then $\ell(c) \leq \lambda n$. 
Let $H : [0,1]^2 \to X$ be a based null-homotopy of $c$ with filling length at most $1+ \FL_X(\lambda n)$.  
By uniform continuity, there exists $\varepsilon >0$ such that $\varepsilon^{-1} \in \mathbb{Z}$ 
and $d_X(H(a),H(b)) \leq 1$ for all $a,b \in [0,1]^2$ with $d_{\mathbb{E}^2}(a,b) \leq \varepsilon$.

Subdivide $[0,1]^2$ into $\varepsilon^{-1}$ rectangles separated by the lines $t = t_j$ where $t_j = j \varepsilon$ and $j=0, 1, \ldots,  \varepsilon^{-1}$.  For all such  $t =t_j$, 
take $0  = s_{t,0} < s_{t,1}  < \ldots < s_{t,k_t} = 1$ in such a way that for all $i$, the restriction of $H_t$ to $[s_{t,i}, s_{t,i+1}]$ is an arc of length at most $\lambda$ and $k_t \leq 1 + \ell(H_t)/\lambda$.  
Mark the points $s_{t,0}, \ldots, s_{t,k_t}$ on each of the lines $t=t_j$.  
Then, for all $j=0, 1, \ldots, \varepsilon^{-1} -1$ and all $i = 1,2, \ldots, k_{t_j} - 1$, join $(s_{t,i},t_j)$ to $(s_{t,i'},t_{j+1})$ by a straight-line segment where $(s_{t,i'},t_{j+1})$ is the first marked point reached from  $(s_{t,i}, t_{j+1})$ by increasing the $s$-coordinate.  
Note that 
\begin{equation} \label{go right}
d( H_{t_j}( s_{t,i}), H_{t_{j+1}}( s_{t,i'}))  \ \leq \  1 + \lambda.
\end{equation}
In the same way, for all $j=1, \ldots, \varepsilon^{ -1}$ and  for all $i = 1,\ldots, k_{t_j} - 1$ such that $(s_{t,i},t_j)$ is not the terminal vertex of one of the edges we just connected, join $(s_{t,i},t_j)$ to some $(x_{t,i'},t_{j-1})$ for which 
\begin{equation} \label{go left} 
d(H_{t_j}( s_{t,i}), H_{t_{j-1}}( s_{t,i'}))  \ \leq \  1 + \lambda,
\end{equation}
so as to produce a diagram $\Delta$ in which every 2-cell has boundary circuit of combinatorial length at most $4$.  
  
Orient every edge  $e$ of $\Delta$ arbitrarily and define $$g_e \ := \ \Psi(H_{t'}(s')) \Psi(H_t(s))^{-1}$$ when the initial and terminal points of $e$ are $(s,t)$ and $(s',t')$, respectively.  It follows from (\ref{go right}) and (\ref{go left}) that in the word metric associated to $\AA$ the distance from $1$
to $g_e$, denoted $\abs{g_e}$, is at most a constant $K=K(\AA)$. Subdivide $e$ into a path of $\abs{g_e}$ edges; give each of these new edges an  orientation and a labelling by a letter in $\AA$ so that one reads a word representing $g_e$ along the path.   
Make all the choices in the construction above in such a way that $w$ labels the line $t=0$ and all the other edges in  $\partial \Delta$ are labelled by $e$. 

The shelling $(\Delta_i)$ of $\Delta$ which strips away the rectangles from left to right, shelling each in turn from top to bottom, has $$\max_i \ell(\Delta_i) \ \leq \  K(1+ \FL_X(\lambda n)) + 4K.$$
Let $w_i$ be the word one reads around the boundary circuit of $\Delta_i$.
Each 2-cell in each $\Delta_i$ has boundary circuit labelled by a null-homotopic word that may not be in $\RR$, but has length at most $4K$.  So it is possible to interpolate between the $w_i$ to produce a null-sequence (as in Section~\ref{filling length}) for $w = w_0$ with respect  to $\PP$; in this sequence  every word has length at most   $K(1+ \FL_X(\lambda n)) + 4K$ plus a universal constant.
Thus $\FL_{\PP}  \ \preceq  \ \FL_X$, as required. 

That $\FFL_{\PP} \preceq \FFL_X$ can be proved in the same way.  The argument needs to be developed further to show that $\FFFL_{\PP} \preceq \FFFL_X$.  Given the word-like loop $c$, one takes a free and fragmenting null-homotopy $H : S \to X$ of $c$ with $\ell(H)$ at most $\FFFL_X(\ell(c))+1$ \emph{and} with the property that whenever loops of length less than some prior fixed constant  appear, those loops are contracted to points before any further bifurcations occur.      
This implies that for all $t$,   the number of connected components  $k_t$ in the fibre $S_t= \bigsqcup_{i=1}^{k_t} I_{t,i} \times \set{t}$ is at most a constant times $(1 + \ell(H_t))$.  For the construction of $\Delta$ we inscribe $S$ with the arcs of its intersection with the lines  $t=t_j$ and with the additional  lines $t=\tau_l$, where $\tau_l$ are the critical $t$-values of $H$.  Using $H_t$ and $\Psi$ as before, we subdivide the fibre $S_t$ into edges and label each of its $k_t$ connected components by a null-homotopic word  -- this works as before, except we additionally insist that the end points of the closed intervals $I_{t,i}$ comprising $S_t$ be included amongst the $s_i$ -- this may add  $k_t$ to the total length the words along $S_t$, but the argument given above ensures that this additional cost is no more than a constant times $(1+ \ell(H_t))$.             
\end{proof}

\bibliographystyle{plain}
\bibliography{bibli}

\small{ \ni \textsc{Martin R.\ Bridson} \rule{0mm}{6mm} \\
Mathematics, Huxley Building, Imperial College, London, SW7 2AZ,
UK \\ \texttt{m.bridson@imperial.ac.uk, \
http:/\!/www.ma.ic.ac.uk/$\sim$mbrids/ }

\ni  \textsc{Tim R.\ Riley} \rule{0mm}{6mm} \\ Mathematics
Department, 310 Malott Hall, Cornell University, Ithaca, NY
14853-4201, USA \\ \texttt{tim.riley@math.cornell.edu, \
http:/\!/www.math.cornell.edu/$\sim$riley/ } }\end{document}